\documentclass{amsproc}

\usepackage{verbatim} 
\usepackage{amsmath}
\usepackage{amsfonts}
\usepackage{amssymb}
\usepackage{mathrsfs}
\usepackage{amsthm}
\usepackage{newlfont}

\usepackage{fullpage,amsmath, amssymb,amscd}
\usepackage{latexsym, epsfig, color}
\usepackage[mathscr]{eucal}
\usepackage{xypic, graphicx}
\usepackage{amscd}
\usepackage[all, cmtip]{xy}

\usepackage{float, hyperref, nccmath, multirow, caption, subcaption}

\newcommand\cyr{%
 \renewcommand\rmdefault{wncyr}%
 \renewcommand\sfdefault{wncyss}%
 \renewcommand\encodingdefault{OT2}%
\normalfont\selectfont} \DeclareTextFontCommand{\textcyr}{\cyr}

\newtheorem{theorem}{Theorem}

\def\Z{\mathbb Z}

\def\Q{\mathbb Q}

\def\F{\mathbb F}

\def\End{\operatorname{End}}

\def\Gal{\operatorname{Gal}}

\def\mod{\operatorname{mod}}

\def\GL{\operatorname{GL}}

\def\PGL{\operatorname{PGL}}

\def\O{\operatorname{O}}
\def\o{\operatorname{o}}

\def\log{\operatorname{log}}

\def\ds{\displaystyle}

\begin{document}

\title{
The absolute discriminant of the endomorphism ring of
most reductions of a non-CM elliptic curve is close to maximal 
}

\date{February 25,  2020}

\author{
Alina Carmen Cojocaru and Matthew Fitzpatrick
}
\address[Alina Carmen  Cojocaru]{
\begin{itemize}
\item[-]
Department of Mathematics, Statistics and Computer Science, University of Illinois at Chicago, 851 S Morgan St, 322
SEO, Chicago, 60607, IL, USA;
\item[-]
Institute of Mathematics  ``Simion Stoilow'' of the Romanian Academy, 21 Calea Grivitei St, Bucharest, 010702,
Sector 1, Romania
\end{itemize}
} \email[Alina Carmen  Cojocaru]{cojocaru@uic.edu}

\address[Matthew Fitzpatrick]{
\begin{itemize}
\item[-]
Department of Mathematics, Statistics and Computer Science, University of Illinois at Chicago, 851 S Morgan St, 322
SEO, Chicago, 60607, IL, USA;
\end{itemize}
} \email[Matthew Fitzpatrick]{mfitzp6@uic.edu}

\renewcommand{\thefootnote}{\fnsymbol{footnote}}
\renewcommand{\thefootnote}{\arabic{footnote}}

\renewcommand{\thefootnote}{\fnsymbol{footnote}}
\renewcommand{\thefootnote}{\arabic{footnote}}

\thanks{A.C.C. was partially supported  by  a Collaboration Grant for Mathematicians from the Simons Foundation  
under Award No. 318454.}

\thanks{February 25, 2020}

\begin{abstract}
Let $E/\Q$ be a non-CM elliptic curve. Assuming GRH, we prove that, for a 
set of primes $p$ of density $1$, the absolute discriminant of the $\F_p$-endomorphism ring of the 
reduction  of $E$ modulo $p$ is close to maximal.
\end{abstract}

\maketitle

\section{Introduction}

Let $E/\Q$ be an elliptic curve defined over the field of rational numbers, of conductor $N_E$, and let
$p \nmid N_E$ be a prime of good reduction for $E$. We denote by $E_p/\F_p$ the reduction of $E$ modulo $p$
and we recall that it is an elliptic curve defined over the finite field $\F_p$ with $p$ elements, with the property that
$|E_p(\F_p)| = p + 1 - a_p$ for some integer $a_p$ satisfying $|a_p| < 2 \sqrt{p}$. Consequently, the polynomial
$X^2 - a_p X + p$ has two complex conjugate roots, $\pi_p$ and $\overline{\pi}_p$, satisfying
$|\pi_p| = \sqrt{p}$. 
Upon identifying any one of these roots, say $\pi_p$,  with the $p$-th power Frobenius endomorphism of $E_p/\F_p$,
we obtain  the embeddings of  imaginary quadratic orders
$\Z[\pi_p] \leq \End_{\F_p}(E_p) \leq {\cal{O}}_{\Q(\pi_p)}$
in  the field $\Q(\pi_p)$, with ${\cal{O}}_{\Q(\pi_p)}$ denoting the maximal order of $\Q(\pi_p)$.
Focusing on the discriminants of these orders, we obtain the relation
\begin{equation}\label{orders-disc}
a_p^2 - 4 p = b_p^2 \Delta_p,
\end{equation}
where $\Delta_p$ denotes the discriminant of  $\End_{\F_p}(E_p)$ and $b_p$ denotes the unique positive integer
satisfying the $\Z$-module isomorphism $\End_{\F_p}(E_p)/\Z[\pi_p] \simeq \Z/b_p \Z$.
If $p \geq 5$ is supersingular, then $\Delta_p \in\{ -p, -4 p\}$,
while if $p$ is ordinary and $\End_{\overline{\Q}}(E) \not\simeq \Z$, then $\Delta_p$ equals the
discriminant of the imaginary quadratic order $\End_{\overline{\Q}}(E)$. 

The goal of this article is
 to focus on  the setting $p$ ordinary and $\End_{\overline{\Q}}(E) \simeq \Z$
and
to investigate the growth of the absolute discriminant $|\Delta_p|$ as a function of $p$, in
particular in relation to the upper bound $4 p - a_p^2$ arising from (\ref{orders-disc}).
In this setting, it was shown in \cite{Sc89} that
$|\Delta_p|$ does indeed grow  with $p$: there exists a positive constant $c(E)$ such that, for any prime $p \nmid N_E$,
\begin{equation*}\label{schoof-1}
|\Delta_p| \geq c(E) \frac{(\log p)^2}{(\log \log p)^4}.
\end{equation*}
Moreover, it was shown in \cite{Sc89}  
under the assumption of the Generalized Riemann Hypothesis (GRH for short)
that
there exists a positive constant $c'(E)$ 
and
there exist infinitely many primes $p$ 
such that
\begin{equation*}\label{schoof-2}
|\Delta_p| \leq c'(E) p^{\frac{2}{3}} \log p.
\end{equation*}
Under similar hypotheses, we will prove that, in fact,
 the growth of $|\Delta_p|$ is very close to the growth of $4p - a_p^2$
 for most primes:

\begin{theorem}\label{main-thm}
Let $E/\Q$ be an elliptic curve of conductor $N_E$ with $\End_{\overline{\Q}}(E) \simeq \Z$.
Assume that GRH holds for the division fields of $E$.
Then, for any function $f: (0, \infty) \longrightarrow (0, \infty)$ 
with $\ds\lim_{x \rightarrow \infty} f(x) = \infty$, 
\begin{equation}\label{main-asymp}
\#\left\{
p \leq x: p \nmid N_E, 
|\Delta_p| \geq \frac{4 p - a_p^2}{f(p)}
\right\}
\sim
\pi(x),
\end{equation}
where $\pi(x)$ denotes the number of primes up to $x$.
\end{theorem}

The growth of $|\Delta_p|$ has also been investigated in other settings, 
including
 that of  arbitrary elliptic curves over finite fields 
 -- see \cite{LuSh07}, \cite{Sh10}, \cite{Sh15} --
 and that  of finite Drinfeld modules -- see \cite{CoPa20}.

Regarding Theorem \ref{main-thm},
the proximity of $|a_p|$ to $\left[2 \sqrt{p}\right]$  was studied in several papers by K. James and
his co-authors, such as  \cite{JaTrTrWeZa16} and \cite{GiJa18} (see also the recent follow-up \cite{DaGaMaPrTB20}).
In \cite{JaTrTrWeZa16}, it is conjectured that, when $\End_{\overline{\Q}}(E) \simeq \Z$,
the number of primes $p \leq x$ with $|a_p| = \left[2 \sqrt{p}\right]$,
called extremal primes, is asymptotically equal to $C(E) \frac{x^{\frac{1}{4}}}{\log x}$ for some constant $C(E)$;
in \cite{GiJa18}, it is proved that this conjecture holds on average over two-parameter families of elliptic curves $E/\Q$
(the majority of which have a trivial endomorphism ring $\End_{\overline{\Q}}(E)$). 
Thus extremal primes are not expected to contribute to the left hand side of  (\ref{main-asymp}).

The proof of Theorem \ref{main-thm} relies on the intimate connection between the integer $b_p$
and the discriminant $\Delta_p$ provided by equation (\ref{orders-disc}), as well as on a characterization criterion of the divisors of $b_p$ through splitting conditions on $p$ in certain subfields of the division fields of $E$. Thanks to these connections, we approach the study of the growth of $|\Delta_p|$ as a potential application of the Chebotarev Density Theorem in an infinite family of number fields.  As such, the assumption of  GRH facilitates best possible error terms.
Even under this assumption, the accumulation of all occurring error terms is overbearing. This we circumvent by resorting to an application of the Square Sieve, which, itself, incorporates another application of the Chebotarev Density Theorem.

{\bf{Notation.}}
In what follows, we use the standard $\o$, $\O$, $\ll$, and $\asymp$ notation:
given suitably defined real functions $h_1, h_2$,
we say that 
$h_1 = \o(h_2)$ if $\ds\lim_{x \rightarrow \infty} \frac{h_1(x)}{h_2(x)} = 0$;
we say that
$h_1 = \O(h_2)$ or $h_1 \ll h_2$ 
if 
$h_2$ is positive valued 
and
 there exists a positive constant $C$ such that 
$|h_1(x)| \leq C h_2(x)$ for all $x$ in the domain of $h_1$;
we say that
$h_1 \asymp h_2$ 
if 
$h_1$, $h_2$ are positive valued 
and 
$h_1 \ll h_2 \ll h_1$;
we say that
$h_1 = \O_D(h_2)$ or $h_1 \ll_D h_2$
if 
$h_1 = \O(h_2)$
and
the implied $\O$-constant $C$ depends on priorly given data  $D$;
similarly,
we say that 
$h_1 \asymp_D h_2$ 
if the implied constant in either one of the $\ll$-bounds
$h_1 \ll h_2 \ll h_1$
depends on priorly  given data $D$.

\section{Proof of the main theorem}

Let $E/\Q$ be an elliptic curve of conductor $N_E$, with $\End_{\overline{\Q}}(E) \simeq \Z$.
Let $f: (0, \infty) \longrightarrow (0, \infty)$ be a function satisfying 
$\ds\lim_{x \rightarrow \infty} f(x) = \infty$.
Without loss of generality, we may assume that 
$f(x)$ is increasing, for we may replace $f(x)$ with $\sup\{f(z):z\le x\}$.
With notation as in Section 1, we observe that, thanks to (\ref{orders-disc}), in order to prove (\ref{main-asymp})
it is enough to prove 
\begin{equation}\label{main-asymp-o}
\#\left\{
p \leq x: p \nmid N_E, b_p > \sqrt{f(p)}
\right\}
= \o (\pi(x)).
\end{equation}
This we will do by exploring the divisibility properties of $b_p$.

As usual, for a positive integer $n$, we denote by $E[n]$ the group of $\overline{\Q}$-rational points of $E$
of order dividing $n$ and by $\Q(E[n])$ the finite, Galois extension of $\Q$ generated by the $x$ and $y$ coordinates
of the points of $E[n]$. We view the Galois group $\Gal(\Q(E[n])/\Q)$ as a subgroup of $\GL_2(\Z/n \Z)$
under the residual modulo $n$ Galois representation of $E$. 
With this notation, the main result of \cite{DuTo02} specialized to  elliptic curves over $\Q$
 states that, for any prime $p \nmid n N_E$,  
the reduction modulo $n$ of
the integral matrix
$$
\left(
\begin{matrix}
\frac{a_p + b_p \delta_p}{2} & b_p
\\
\frac{b_p (\Delta_p - \delta_p)}{4} & \frac{a_p - b_p \delta_p}{2}
\end{matrix}
\right),
$$
with $\delta_p := 0, 1$ according to whether $\Delta_p \equiv 0, 1 (\mod 4)$,
gives a 
representative of the conjugacy class defined by the Artin symbol
$\left(\frac{\Q(E[n])/\Q}{p}\right)$ in $\Gal(\Q(E[n])/\Q)$.
Consequently, 
upon defining
$$
J_n :=
\left\{z \in \Q(E[n]): \sigma(z) = z \quad \forall \sigma \in \Gal(\Q(E[n])/\Q) \ \text{a scalar element}
\right\},
$$
we obtain the criterion
\begin{equation}\label{n-divides-b_p-splits}
n \mid b_p 
\
\Leftrightarrow
\
p \ \text{splits completely in} \ J_n/\Q.
\end{equation}

For each prime $p$, there are unique positive integers $r_p$ and $m_p$, with $m_p$ squarefree, such that 
$$
b_p^2|\Delta_p|=r_p^2m_p.
$$ 
Observe that we must have $b_p \mid r_p$, which gives $b_p \leq r_p$. 
Recalling (\ref{orders-disc}),  observe that
\begin{equation}\label{n-divides-b_p-square}
4 p - a_p^2 = r_p^2 m_p,
\end{equation}
which gives
$r_p < 2 \sqrt{p}$
and
\begin{equation}\label{m-bound}
\ds m_p=\frac{4p-a_p^2}{r_p^2}\leq\frac{4p}{b_p^2}.
\end{equation}
Furthermore, observe that the divisibility $n \mid b_p$ implies that $n \mid r_p$ and, in particular, that $n\leq r_p$.

Now let us proceed to bounding from above the left hand side of (\ref{main-asymp-o}).
We fix an arbitrary parameter $z = z(x)$ satisfying $0 < z < x$ and  define
$$
g(z) := \inf\left\{f(p): z< p \leq x\right\}.
$$
Note that $f(z) \leq g(z)$.
We have the bounds
\begin{eqnarray}\label{sum-of-S}
\#\left\{
p \leq x:
p \nmid N_E, b_p > \sqrt{f(p)}
\right\}
&\leq&
\pi(z)+
\#\left\{
z<p \leq x:
p \nmid N_E, b_p > \sqrt{f(p)}
\right\}
\nonumber
\\
&\leq&
\pi(z)+
\#\left\{
z<p \leq x:
p \nmid N_E, b_p > \sqrt{g(z)}
\right\}
\nonumber
\\
&\leq&
\pi(z)
+
\ds\sum_{\sqrt{g(z)} < n \leq 2 \sqrt{x}}
\#\left\{
p \leq x:
p \nmid N_E, b_p = n
\right\}
\nonumber
\\
&\leq&
\pi(z)
+
S\left(\sqrt{g(z)}, 2 \sqrt{x}\right),
\end{eqnarray}
where, for any $u < v \leq 2 \sqrt{x}$,
$$
S(u, v) := \ds\sum_{u < n \leq v}
\#\left\{
p \leq x:
p \nmid N_E, b_p \equiv 0 (\mod  n)
\right\}.
$$

To estimate $S\left(\sqrt{g(z)}, 2 \sqrt{x}\right)$,
it is natural to use criterion (\ref{n-divides-b_p-splits}) and employ an effective version
of the Chebotarev Density Theorem as proved in \cite{LaOd77}.
However,  this strategy leads to limitations on the range of $n$
-- indeed, looking at (\ref{chebotarev-applic}) below,
we see that
an application of the effective version of the Chebotarev Density Theorem to the entire sum
$S\left(\sqrt{g(z)}, 2 \sqrt{x}\right)$  gives rise to an error term of $\O_E(x \log x)$.
As such, we adjust our strategy as follows:
 we fix an arbitrary parameter $y = y(x)$ with $\sqrt{g(z)} \leq y \leq 2 \sqrt{x}$
and analyze each of the sums $S\left(\sqrt{g(z)}, y\right)$ and $S\left(y, 2 \sqrt{x}\right)$ separately,
using criterion (\ref{n-divides-b_p-splits}) for the first sum
and using observation  (\ref{n-divides-b_p-square})
for the second sum.

Specifically,
rewriting $S\left(\sqrt{g(z)}, y\right)$ via (\ref{n-divides-b_p-splits})
and using the effective Chebotarev Density Theorem, conditional upon GRH,
 as  in \cite[Prop. 4.2, p. 523]{CoDu04}, we obtain that
\begin{eqnarray}\label{chebotarev-applic}
S\left(\sqrt{g(z)}, y\right)
&=&
\ds\sum_{\sqrt{g(z)} < n \leq y}
\frac{1}{[J_n : \Q]} \pi(x)
+
\O\left(
\ds\sum_{\sqrt{g(z)} < n \leq y} 
x^{\frac{1}{2}} \log (N_E n x)
\right)
\nonumber
\\
&=&
\ds\sum_{\sqrt{g(z)} < n \leq y}
\frac{1}{[J_n : \Q]} \pi(x)
+
\O_E\left(
y x^{\frac{1}{2}} \log  (y x)
\right)
\nonumber
\\
&=&
\ds\sum_{\sqrt{g(z)} < n \leq y}
\frac{1}{[J_n : \Q]} \pi(x)
+
\O_E\left(
y x^{\frac{1}{2}} \log  x
\right),
\end{eqnarray}
where in the last line we also used that $y < 2 \sqrt{x}$.

We estimate  the first term on the right hand side of the above equation
by using that,
thanks to Serre's Open Image Theorem \cite{Se72},
the assumption $\End_{\overline{\Q}}(E) \simeq \Z$
gives 
\begin{equation}\label{degree-J_n}
n^3 \ll_E [J_n: \Q] \ll n^3.
\end{equation}
The upper bound is an immediate consequence of the inequality $[J_n:\Q] \leq \left|\PGL_2(\Z/n \Z)\right|$,
while the lower bound requires a more detailed analysis, which is presented in \cite{CoJo20}.

Together with Chebyshev's Theorem, 
the lower bound in (\ref{degree-J_n}) gives that
\begin{eqnarray}\label{S-one}
S\left(\sqrt{g(z)}, y\right)
&\ll_{E}&
\ds\sum_{\sqrt{g(z)} < n \leq y}
\frac{\pi(x)}{n^3} 
+
y x^{\frac{1}{2}} \log x
\nonumber
\\
&\ll_{E}&
\ds\pi(x)\int_{\sqrt{g(z)}}^{y+1}
\frac{1}{t^3}\, dt
+
y x^{\frac{1}{2}} \log x
\nonumber
\\
&\ll_{E}&
\frac{x}{g(z)\log x} 
+
y x^{\frac{1}{2}} \log  x.
\end{eqnarray}

To estimate $S\left(y, 2 \sqrt{x}\right)$,
we use (\ref{n-divides-b_p-square}) - (\ref{m-bound})
and obtain that
\begin{eqnarray*}
S\left(y, 2 \sqrt{x}\right)
&=&
\ds\sum_{y < n \leq 2 \sqrt{x}}
\#\left\{
p \leq x:
p \nmid N_E, b_p \equiv 0 (\mod  n)
\right\}
\\
&\leq&
\ds\sum_{y < n \leq 2 \sqrt{x}}
\#\left\{
p \leq x:
p \nmid N_E, r_p \equiv 0 (\mod  n)
\right\}
\\
&=&
\ds\sum_{y < n \leq 2 \sqrt{x}}\sum_{\substack{n < r \leq 2 \sqrt{x}\\ n\mid r}}
\#\left\{
p \leq x:
p \nmid N_E, r_p=r
\right\}
\\
&\leq&
\ds\sum_{y < r \leq 2 \sqrt{x}}\tau(r)
\#\left\{
p \leq x:
p \nmid N_E, r_p=r
\right\}
\\
&=&
\ds\sum_{y < r \leq 2 \sqrt{x}}
\tau(r)
\#\left\{p \leq x:
p \nmid N_E, 4 p - a_p^2 = r^2 m_p
\right\}
\\
&\leq&
 \ds\sum_{\substack{m \leq \frac{4x}{y^2}\\m \ \text{squarefree}}}\sum_{y < r \leq 2 \sqrt{x}}\tau(r)
\#\left\{p \leq x:
p \nmid N_E, 4 p - a_p^2 = r^2 m
\right\}
\\
&=&
 \ds\sum_{\substack{m \leq \frac{4x}{y^2}\\m \ \text{squarefree}}}\sum_{y < r \leq 2 \sqrt{x}}\tau(r)
\#\left\{p \leq x:
p \nmid N_E, m(4 p - a_p^2) = (mr)^2
\right\},
\end{eqnarray*}
where $\tau(\cdot)$ is the divisor function.
Using the standard bound $\tau(r) \ll_{\varepsilon} r^{\varepsilon}$ for any $\varepsilon >0$
and observing that, 
for each squarefree $m$ and each prime $p$, 
there is at most one $y<r\le 2\sqrt{x}$ such that $4p-a_p^2=r^2m$,
we deduce that
\begin{equation}\label{S-two-one}
S\left(y, 2 \sqrt{x}\right)
\ll_{\varepsilon}
x^{\varepsilon}
\ds\sum_{
m \leq \frac{4 x}{y^2}
\atop{
m \ \text{squarefree}
}
}
S_m(x),
\end{equation}
where
$$
S_m(x) :=
\#\left\{
p \leq x: p \nmid N_E, m \left(4 p - a_p^2\right) \ \text{is a square}
\right\}.
$$

The sum of the counting functions $S_m(x)$ can be estimated 
as an indirect application  of the conditional effective
Chebotarev Density Theorem of \cite{LaOd77}
via the Square Sieve of \cite{HB84} (see also the version of this sieve in \cite{CoFoMu05}),
by relying  once again on GRH and  on the assumption $\End_{\overline{\Q}}(E) \simeq \Z$.
Such an application was pursued in \cite{CoDu04}, \cite{CoFoMu05}, \cite{CoDa08},
and \cite{Sh10}.
While any of the estimates proved in these papers suffices to conclude the proof of our theorem,
we use the best among all of them, which is the one of  \cite[Lem. 3.2, p. 260]{Sh10};
its proof relies on  \cite[Thm. 3, p. 4]{CoDa08} and \cite[Thm. 1, p. 237]{HB95}.
 We deduce that
\begin{equation}\label{S-two-two}
\ds\sum_{
m \leq \frac{4 x}{y^2}
\atop{
m \ \text{squarefree}
}
}
S_m(x)
\ll_{\varepsilon}
\frac{
x^{\frac{13}{7} + \varepsilon}
}{
y^{\frac{13}{7}}
}.
\end{equation}

Next, by putting together
(\ref{sum-of-S}),
(\ref{S-one}),
(\ref{S-two-one}),
and
(\ref{S-two-two}),
we obtain that
\begin{eqnarray*}\label{sum-of-S-revisited}
\#\left\{
p \leq x:
p \nmid N_E, b_p > \sqrt{f(p)}
\right\}
&\ll_{E, \varepsilon}&
\frac{z}{\log z}
+
\frac{x}{g(z)\log x} 
+
y x^{\frac{1}{2}} \log x
+
\frac{
x^{\frac{13}{7} + \varepsilon}
}{
y^{\frac{13}{7}}
}.
\end{eqnarray*}
The trivial endomorphism
 assumption is reflected in each of 
 the second, third and fourth terms
on the right hand side of the above estimate,
while the GRH assumption is reflected  in each of the third and fourth terms.

Finally, we  choose
$y \asymp x^{\frac{19}{40}}$
and
$z \asymp x^\frac{1}{2}$,
we recall that $f(z) \leq g(z)$,
and 
we use the assumption that $f$ is increasing. We deduce that
\begin{eqnarray*}
\#\left\{
p \leq x:
p \nmid N_E, b_p > \sqrt{f(p)}
\right\}
&\ll_{E, \varepsilon}&
\frac{x^\frac{1}{2}}{\log x}
+
\frac{x}{f\left(x^\frac{1}{2}\right) \log x}
+
x^{\frac{39}{40} + \varepsilon}.
\end{eqnarray*}
Note that, since $\ds\lim_{x \rightarrow \infty} f(x) = \infty$, 
the right hand side of the above estimate is $\o(\pi(x))$.
This completes the proof.

\bigskip

\noindent
{\bf{Acknowledgments.}}
We thank Nathan Jones  for helpful comments on a previous draft of this manuscript.
We thank the anonymous referee for suggestions that added to the clarity of the paper.


{\small{

}}

\end{document}